\documentclass[11pt]{amsart}
\pdfoutput=1
\usepackage{amsmath, amssymb, amsthm, amsbsy, mathtools, bbm, amsfonts, color, verbatim, array, floatflt, yfonts, mathrsfs}

\usepackage[margin=1.1in]{geometry}
\usepackage[all]{xy}
\usepackage[english]{babel}
\usepackage[shortlabels]{enumitem}
\usepackage{ stmaryrd } % for double brackets

\usepackage{graphicx, color, overpic}
\usepackage[usenames,dvipsnames,svgnames,table]{xcolor}
\usepackage[pdfpagelabels]{hyperref}
\usepackage{cleveref}
\usepackage{multirow}
\usepackage{easybmat}
\usepackage{adjustbox}

\usepackage{tikz}

\theoremstyle{plain}
\newtheorem{thm}{Theorem}[section]
\newtheorem*{thm*}{Theorem}
\newtheorem*{prop*}{Proposition}

\newtheorem*{lemma*}{Lemma}

\theoremstyle{definition}

\newtheorem{example}[thm]{Example}

\newtheorem{rmk}[thm]{Remark}

\newcommand{\Z}{\mathbb{Z}}

\newcommand{\bE}{\mathbb{E}} % Euclidean space

\newcommand{\aW}{W} % affine Coxeter group
\newcommand{\sW}{W_0} % spherical Weyl group
\newcommand{\xconj}{[x]} % conjugacy class of x
\newcommand{\cent}{\operatorname{C}} % centralizer
\newcommand{\coconj}[2]{\operatorname{C}(#1,#2)} % co-conjugation set for #1 and #2
\newcommand{\coconjW}[2]{\operatorname{C}_{\sW}(#1,#2)} % spherical coconjugation set for #1 and #2
\newcommand{\tcCoconj}{\operatorname{C}} % translation-compatible part of the coconjugation set
\newcommand{\Mov}{\textsc{Mov}} % move-set
\newcommand{\Mod}{\textsc{Mod}} % mod-set
\newcommand{\Fix}{\textsc{Fix}} % fix-set
\newcommand{\range}{\operatorname{range}} % range
\newcommand{\Sym}{\operatorname{Sym}}

\definecolor{amethyst}{rgb}{0.6, 0.4, 0.8}
\definecolor{americanrose}{rgb}{1.0, 0.01, 0.24}

\title{Visualizing conjugation in affine Coxeter groups}

\author{Amy Herron}
\address{Amy Herron, School of Mathematics, University of Bristol, Fry Building, Woodland Road, Bristol, BS8 1UG, United Kingdom}
\email{amy.herron@bristol.ac.uk}

\author{Anne Thomas}
\address{Anne Thomas, School of Mathematics \& Statistics, Carslaw Building F07,  University of Sydney NSW 2006, Australia}
\email{anne.thomas@sydney.edu.au}

\date{\today}

\begin{document}


\maketitle

\begin{abstract}
    Affine Coxeter groups are fundamental objects in mathematics and crystallography, and the corresponding tesselations have inspired artists throughout centuries. If two group elements are conjugate, they have very similar algebraic and geometric properties. Using recent structural results of Mili\'cevi\'c, Schwer and the second author, we develop an app to visualize conjugation in affine Coxeter groups in dimensions 2 and 3. We explain the mathematics underlying this visualization, and connect its visual features to the arts.
\end{abstract}

\section{Introduction}

We describe an app we have developed to visualize important structures in affine Coxeter groups. These collections of symmetries are also known as affine Weyl groups and Euclidean reflection groups. Affine Coxeter groups have long been studied in algebra and geometry, and they are a key family in crystallography: in dimensions 2 and 3, every crystallographic group has finite index in some affine Coxeter group (see, for example,~\cite[Appendix~A]{MST-affine}). Two-dimensional or planar crystallographic groups are also called wallpaper groups (see~\cite{Conway}). 
Affine Coxeter groups thus contain all possible symmetries of the planar tesselations which are significant to many art forms. 

\begin{figure}[!ht]
    \centering
    \includegraphics[width=0.8\textwidth]{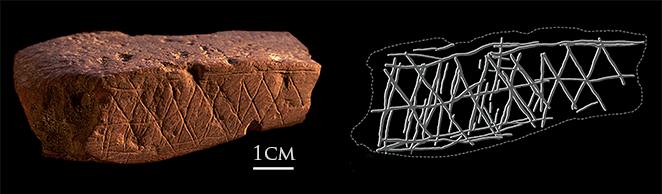}
    \caption{\footnotesize{Ochre stone found at the Blombos Cave site. Author: Chris S. Henshilwood~\cite{Henshilwood-photo}. Licensed under Creative Commons \href{https://creativecommons.org/licenses/by-sa/4.0/deed.en}{Attribution-Share Alike 4.0 International license}.}}
    \label{fig:Blombos}
\end{figure}

Possibly the oldest abstract representation by humans is the engraving from the Blombos Cave in South Africa shown in Figure~\ref{fig:Blombos}, which has been dated to around 70,000 years~\cite{Henshilwood}. 
In visual art, Islamic tilings, such as those in the Alhambra Palace in Granada, Spain  (see~\cite{Bodner} and its references), are key examples of  planar tesselations. Islamic geometry inspired the work of M.C. Escher and Monir Shahroudy Farmanfarmaian (see, for example,~\cite{Shaw}), and continues to inspire contemporary artists, such as Chris Berry~\cite{Berry}.  In fact, the mathematician Coxeter himself directly influenced M.C. Escher's work, particularly Escher's \emph{Circle Limit} series  (see, for example, Coxeter's account in~\cite{Coxeter}). In textile arts, many quilting patterns are  planar tesselations, while in music theory, the Tonnetz maps out a particular tesselation, as explained in the next paragraph.

\begin{figure}[!ht]
    \centering
    \begin{minipage}{.4\textwidth}
        \centering
        \includegraphics[width=0.9\textwidth]{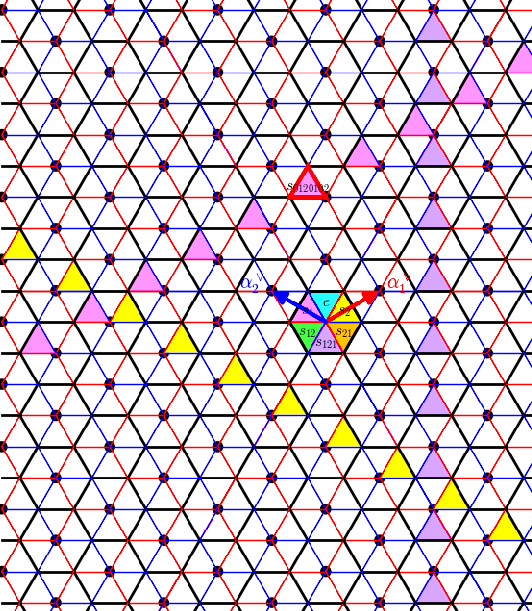}
    \end{minipage}%
    \begin{minipage}{0.3\textwidth}
        \centering
        \includegraphics[width=0.9\textwidth]{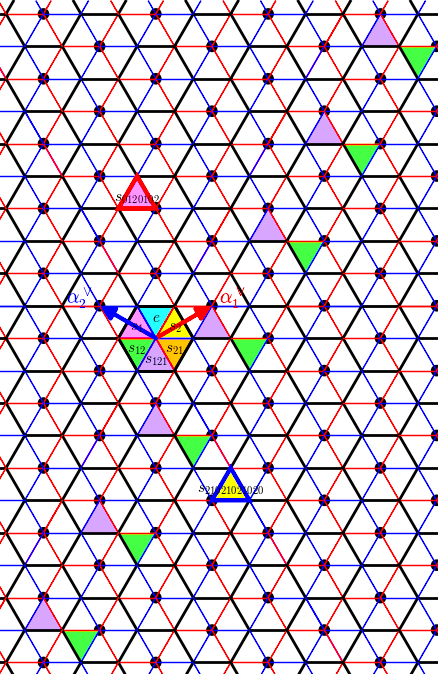}
    \end{minipage}
    \caption{\footnotesize{For $W$ of type $\tilde{A}_2$ and $x = s_{0120102} = t^{2\alpha_1^\vee + 2\alpha_2^\vee}s_1 \in W$, the conjugacy class~$\xconj$ consists of the ``lines" of pink, purple, and yellow tiles depicted on the left. For this same~$x$ and $y = s_{21021021020} = t^{-2\alpha_1^\vee - 3\alpha_2^\vee}s_2 \in \xconj$, the coconjugation set $\coconj x y$ consists of the ``lines" of purple and green tiles shown on the right. The tile~$x$ is outlined in red on both left and right, and $y$ is outlined in blue on the right.}}
    \label{fig:A2tilde_Intro}
\end{figure}

An affine Coxeter group $\aW$ is, by definition, a discrete group of Euclidean isometries generated by reflections in the faces of a convex polytope $P$ in Euclidean space. The action of $\aW$ induces a tiling of Euclidean space by copies of $P$, and the elements of $\aW$ can then be identified with the tiles. A first example is $\aW$ of type $\tilde{A}_2$, which induces the tiling of the Euclidean plane by equilateral triangles, as depicted in  Figure~\ref{fig:A2tilde_Intro} (and perhaps also Figure~\ref{fig:Blombos}).  
This same triangular tiling appears in music theory as the Tonnetz, in which vertices represent pitch classes and triangles represent major and minor triads (see \cite{crans2009musical}).  
Moving a single voice transforms a triad into one of its three neighbors by reflecting the triangle across one of its edges.  Thus, these operations act exactly as reflections in the affine Coxeter group.  Because pitch classes are considered modulo the octave, the Tonnetz is, more precisely, the quotient of this tiling by a finite-index sublattice of the coroot lattice (defined below), so that the planar tiling of Figure \ref{fig:A2tilde_Intro} is its universal cover.

Two elements $x$ and $y$ of an affine Coxeter group $\aW$ are \emph{conjugate} if there is some $z \in \aW$ such that $zxz^{-1} = y$, and the \emph{conjugacy class} of $x \in \aW$ is the set of all conjugates of $x$:
\[
\xconj = \{ zxz^{-1} \mid z \in \aW \}.
\]
Elements of the same conjugacy class have closely related  properties. Algebraically, if $\varphi:\aW \to A$ is a homomorphism from $\aW$ to an abelian group $A$, then for every $y \in [x]$ we have $\varphi(y) = \varphi(x)$. Geometrically, every conjugate of a reflection (respectively, translation) in $\aW$ is another reflection (respectively, translation). The set of elements which conjugate $x$ to itself is called the \emph{centralizer} of $x$:
\[\cent(x) = \{ z \in \aW \mid zxz^{-1} = x\}.\]
More generally, if $x$ and $y$ are conjugate then the corresponding \emph{coconjugation set} is given by
\[\coconj x y = \{ z \in \aW \mid zxz^{-1} = y\}.\]

In this paper we describe an app \texttt{AffineCoxeterExplorer.java} \cite{app}  that we have developed to visualize conjugacy classes and coconjugation sets in affine Coxeter groups in dimensions 2 and 3. Using the identification of the elements of $\aW$ with the tiles in the induced tesselation, we color all tiles which are in the desired conjugacy class or coconjugation set (and which lie inside a bounding box of size chosen by the user). For example, in type $\tilde{A}_2$, Figure~\ref{fig:A2tilde_Intro} shows the output from our app for a conjugacy class $\xconj$ and a coconjugation set $\coconj x y$. (All figures appearing in this paper were generated using version 1.1 of our app.)

Our coloring of certain tiles breaks the overall symmetry of the tesselation. Due to such symmetry-breaking, the use of color in Islamic tilings is sometimes ignored, for example in Bodner's analysis of the planar crystallographic groups appearing in the Alhambra~\cite{Bodner}. Empar\'an~\cite{Emparan} argues, however, that although such color use can break the wallpaper symmetry, it is intentional, and can still maintain balance and some visual symmetry, as well as serving  metaphorical purposes. She writes: ``\dots chromatic asymmetry is neither a
mathematical error nor an oversight \dots chromatic asymmetry would be a symbolic tool to generate a specific
aesthetic experience in the spectator." (see~\cite[p. 303]{Emparan}). 

Similar observations can be made about many quilting patterns, where fabric choices create symmetric and meaningful designs within tesselations. 
By analogy, our use of color makes visible structures that have their own mathematical significance and symmetries.  For instance, the colors in Figure \ref{fig:A2tilde_Intro} are precisely what makes the conjugacy class visible against the uncolored tesselation: the overall tesselation has more symmetry but carries none of the conjugacy class information.

\begin{figure}[!ht]
    \centering
    \begin{minipage}{.3\textwidth}
        \centering
        \includegraphics[width=0.9\textwidth]{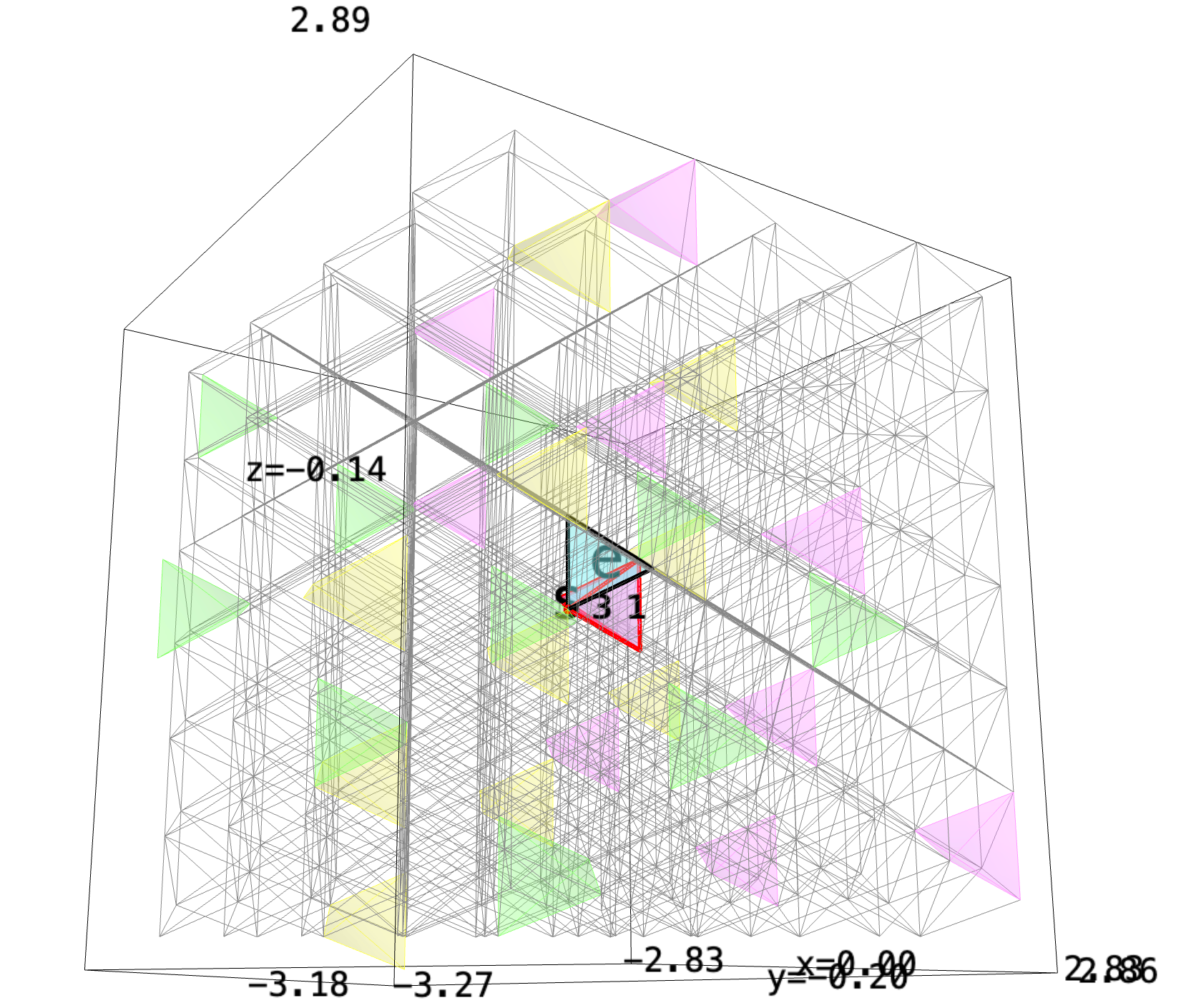}
    \end{minipage}%
    \begin{minipage}{0.3\textwidth}
        \centering
        \includegraphics[width=0.9\textwidth]{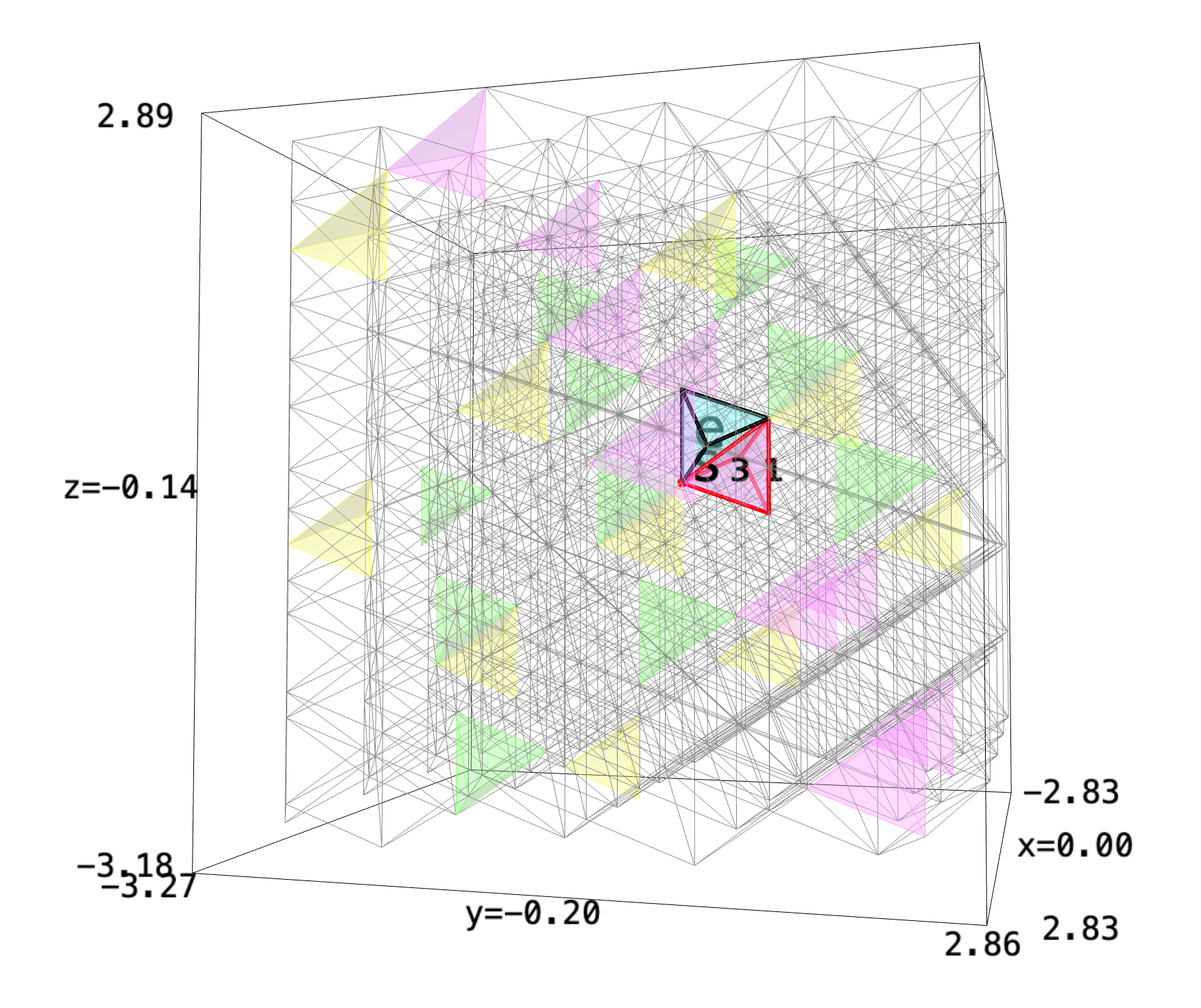}
    \end{minipage}
      \begin{minipage}{0.3\textwidth}
        \centering
        \includegraphics[width=0.9\textwidth]{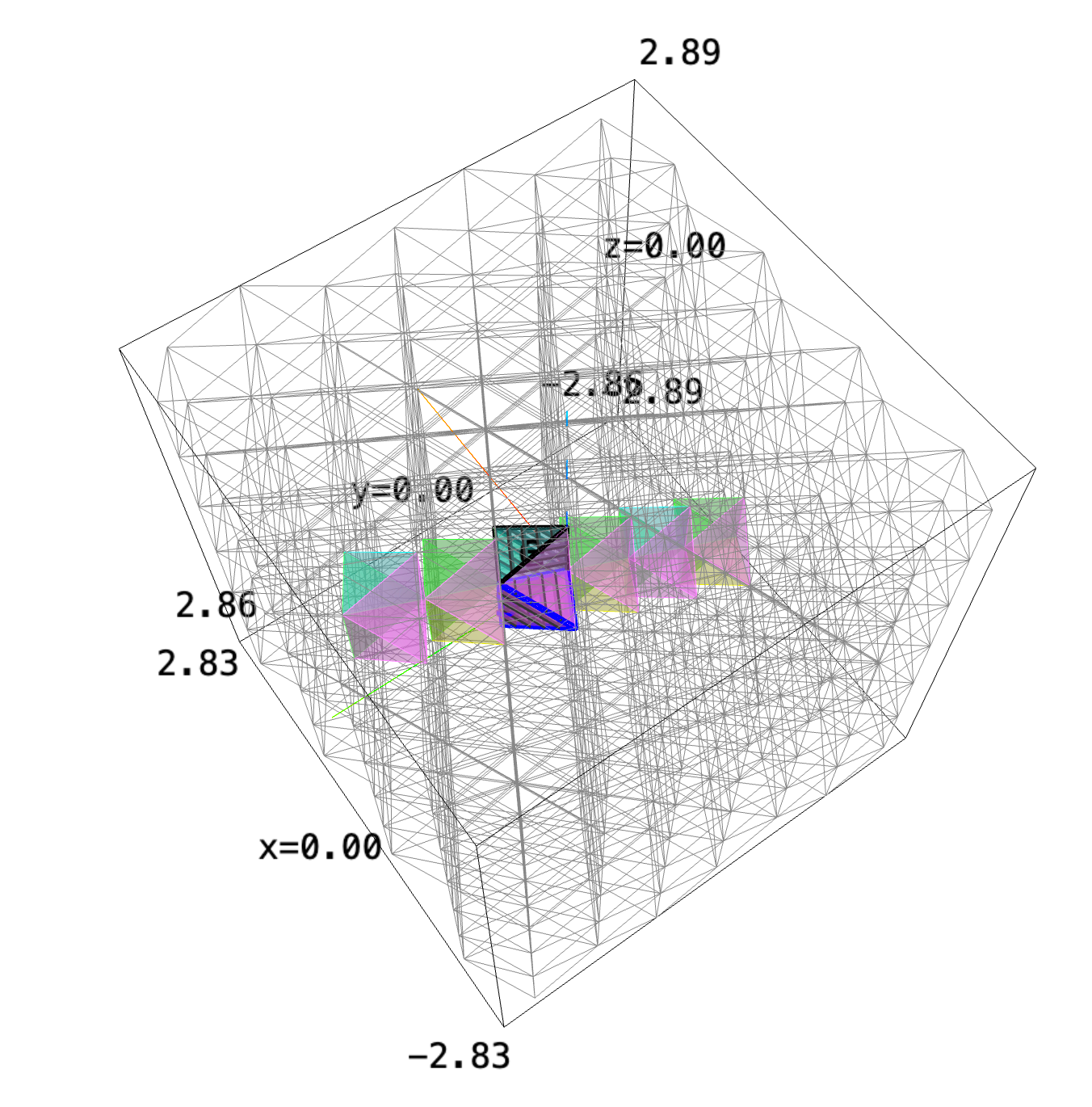}
    \end{minipage}
    \caption{\footnotesize{For $W$ of type $\tilde{A}_3$ and $x = w = s_{13} \in W$, the conjugacy class~$\xconj$ consists of the sheets of pink, green, and yellow tiles depicted in two different views on the left and in the center. The centralizer $\cent(x)$ consists of the lines of pink, green, and blue tiles shown on the right, together with the striped identity alcove (see Section \ref{subsec:features}); here, equal colors indicate conjugate spherical directions
    in $W_0$, and the four colors show that $\cent(x)$ meets exactly four conjugacy classes in~$W_0$. Observe that $\Mov(x)$ is $2$-dimensional and $\Fix(x) = \Fix(w)$ is $1$-dimensional.}}
    \label{fig:A3tilde_Intro}
\end{figure}

Our first goal with this visualization is to illustrate how both conjugacy classes and coconjugation sets lie along certain classically-defined  subspaces. 
For example, in both parts of Figure~\ref{fig:A2tilde_Intro}, the shaded tiles lie along certain lines.   To make this more precise, write~$\bE$ for the Euclidean space tesselated by the action of $\aW$. Then the conjugacy class $\xconj$ lies along certain images of an affine subspace of $\bE$ called the \emph{move-set}:
\[ \Mov(x) = \{ q \in \bE \mid q = xp - p \mbox{ for some }p \in \bE \} = \range(x - I) = (x-I)\bE, \]
where $I$ is the identity map on $\bE$. Each of the ``lines" of tiles on the left of Figure~\ref{fig:A2tilde_Intro} is called a \emph{component} of the conjugacy class, and it can further be seen that there are symmetries between as well as within components. In Figure~\ref{fig:A3tilde_Intro}, the left and center depict the three components of a conjugacy class in type $\tilde{A}_3$.

For coconjugation sets, we note that any $y \in \aW$ can be decomposed uniquely as the product of a translation $t^\lambda$ and an element $w \in \aW$ which fixes the origin. Given $y = t^\lambda w$, the coconjugation set $\coconj x y$ then lies along certain images of a linear subspace of $\bE$ called the \emph{fix-set}:
\[ \Fix(w) = \{ p \in \bE \mid w p = p \} = \ker(w - I). \]
This is illustrated on the right of  Figure~\ref{fig:A2tilde_Intro} and of Figure~\ref{fig:A3tilde_Intro}.

%%%%%%%%%%%%%%%%%%%%%%%%%%%%%%%%%%%%%%%%%%%%%%%%%%
%   Comment on the right-most picture in Figure 2:
%
%   word      permutation      cycle type       color
%   e         (1)               (1,1,1,1)       cyan (striped)
%   s_1       (1,2)             (2,1,1)         pink/magenta
%   s_3       (3,4)             (2,1,1)         pink/magenta
%   s_31      (1,2)(3,4)        (2,2)           green/lime
%   s_2312    (1,3)(2,4)        (2,2)           green/lime
%   s_123121  (1,4)(2,3)        (2,2)           green/lime
%   s_23121   (1,3,2,4)         (4)             blue/cornflower
%   s_12312   (1,4,2,3)         (4)             blue/cornflower
%%%%%%%%%%%%%%%%%%%%%%%%%%%%%%%%%%%%%%%%%%%%%%%%%%

\begin{figure}[!ht]
    \centering
    \includegraphics[width=0.5\textwidth]{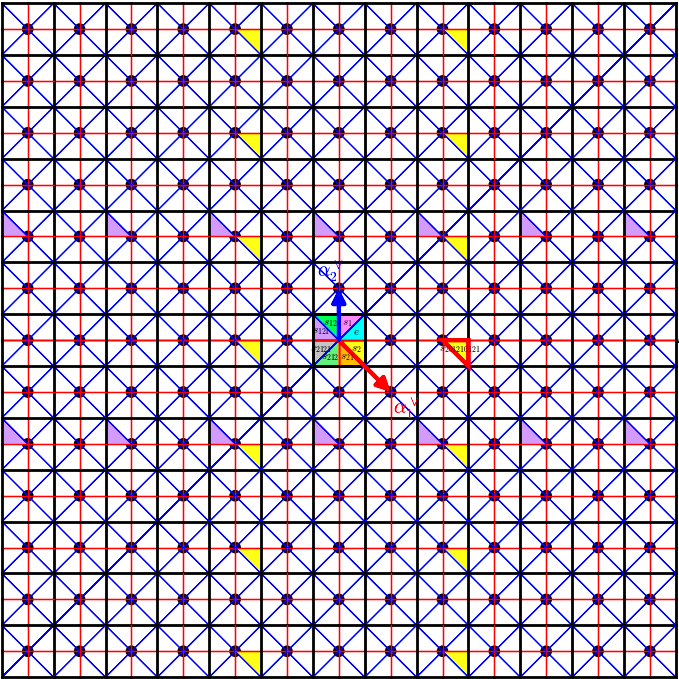}
    \caption{\footnotesize{The conjugacy class of the element $x = s_{201210121} = t^{2\alpha_1^\vee + 2\alpha_2^\vee}s_2$ in~$W$ of type~$\tilde{C}_2$ consists of the horizontal lines of purple tiles, and the vertical lines of yellow tiles. Notice the ``gaps" between shaded alcoves here.}}
    \label{fig:C2tilde_ConjClass_Intro}
\end{figure}

Our second goal is to illustrate some finer periodic behavior  for  conjugacy classes. The \emph{coroot lattice}~$L$ is the set of vectors $\lambda \in \bE$ so that $\aW$ contains the translation by $\lambda$; we indicate the elements of $L$ by heavy black dots in our $2$-dimensional figures. The \emph{mod-set} for an element $x \in \aW$ is then the $L$-analog of the move-set:
\[ \Mod(x) = \{ \mu \in L \mid \mu = x\lambda - \lambda \mbox{ for some }\lambda \in L \} = (x-I)L. \]
Only the tiles along the (images of the) move-set which have translation part in the (images of the) mod-set are actually in the conjugacy class $\xconj$. For example, in Figure~\ref{fig:C2tilde_ConjClass_Intro}, along the horizontal and vertical lines of tiles, only every second translation by an element of $L$ is in~$\xconj$.

The geometric description of conjugacy classes and coconjugation sets that we have just sketched is established in recent work of Mili\'cevi\'c, Schwer and the second author~\cite{MST-affine,MST-Euclidean}. This description actually holds for many other Euclidean isometry groups, including the full isometry group of Euclidean space (a non-discrete group). Thus our app illustrates, for a family of groups important to mathematics, crystallography, the visual arts, and music theory, a much wider connection between algebra and geometry. 

\subsection*{Structure of the paper} We give some additional mathematical background in Section~\ref{sec:background}. We then describe our visualization in Section~\ref{sec:visualization}, including explanations of the choices we made in order to  illustrate mathematical features, and a discussion of our color choices and what they display to the viewer.

\subsection*{Acknowledgements} We thank Chris Berry for discussing his art with us,  Henry Segerman and John Voight for their  advice on publication venues, and Elizabeth Mili\'cevi\'c and Petra Schwer for their comments on this manuscript and our app.    
This work was supported by the Additional Funding Programme for Mathematical Sciences, delivered by EPSRC (EP/V521917/1) and the Heilbronn Institute for Mathematical Research. The second author was supported in part by Australian Research Council Grant FT250100160.

\section{Mathematical background}\label{sec:background}

In this section we give some additional, more precise background. Section~\ref{sec:Coxeter} recalls some of the general theory of  affine Coxeter groups and relates these groups to symmetries of planar tesselations.  Section~\ref{sec:MST} then recalls the relevant results from~\cite{MST-affine,MST-Euclidean}.

\subsection{Affine Coxeter groups}\label{sec:Coxeter}

We recall some key notions concerning affine Coxeter groups. A concise reference for this material is~\cite[Chapter 2]{Ronan}, and the standard reference is~\cite{Humphreys}. As in~\cite{MST-affine,MST-Euclidean}, we follow the conventions of~\cite{Bourbaki}.

A \emph{Coxeter group} is a group $W$ with presentation of the form
\[
W = \langle S \mid (st)^{m_{st}} \mbox{ for all } s,t \in S\rangle
\]
where $m_{ss} = 1$ for all $s \in S$, and $m_{st} = m_{ts} \in \{2,3,4,\dots\} \cup \infty$ if $s, t \in S$ are distinct. Here, $m_{st} = \infty$ means that the product $st$ has infinite order. Given $x \in \aW$, an expression $x = s_{i_1} \dots s_{i_k}$ where $s_{i_1}, \dots, s_{i_k} \in S$ is called a \emph{word} for $x$, and is  a \emph{reduced} word if no shorter product of elements of $S$ equals $x$. We often abbreviate a product $s_{i_1} \dots s_{i_k}$ to $s_{i_1 \dots i_k}$, as in Figures~\ref{fig:A2tilde_Intro}--\ref{fig:C2tilde_ConjClass_Intro}. We write $e$ for the identity element of $W$, which is represented by the empty word.

For $W$ and $S$ as in the previous paragraph, the pair $(W,S)$ is called a \emph{Coxeter system}. 
A Coxeter system $(W,S)$ is \emph{irreducible} if there is no nontrivial partition of $S$ into two disjoint commuting subsets, \emph{spherical} if $W$ is finite, and \emph{affine} if $S$ is the set of reflections in the faces of a compact polytope $P$ in some Euclidean space $\bE$. For convenience, if $\dim(\bE) = n$, we define the corresponding affine Coxeter system $(W,S)$ to be \emph{$n$-dimensional} (this terminology is not standardized). 

The irreducible affine Coxeter systems were classified by Coxeter in the 1930s, and correspond to the case that $P$ is a simplex. Hence if $(W,S)$ is $n$-dimensional and irreducible affine, the generating set $S$ has $n+1$ elements; these are denoted $s_0,s_1,\dots,s_n$. For all affine Coxeter systems, the action of $\aW$ induces a tesselation of $\bE$ by copies of the polytope~$P$. We call these copies of~$P$ the \emph{alcoves}, and identify the elements of $\aW$ with the alcoves in the induced tesselation. 

\begin{example}\label{eg:dimension2}
The $2$-dimensional affine Coxeter systems $(W,S)$ are those where $W$ is generated by the set $S$ of reflections in the sides of a polygon $P$ in the Euclidean plane. We list all such systems below by their classical type, and record the corresponding polygon~$P$. 
\begin{enumerate}
    \item Type $\tilde{A}_1 \times \tilde{A}_1$: $P$ is a square.
    \item Type $\tilde{A}_2$: $P$ is an equilateral triangle, as seen in Figure~\ref{fig:A2tilde_Intro}.
    \item Type $\tilde{B}_2$ and type $\tilde{C}_2$: in both these types, $P$ is a right-angled isosceles triangle. (The coroot lattices in these types are not the same.) See Figure~\ref{fig:C2tilde_ConjClass_Intro} for type $\tilde{C}_2$. 
    \item Type $\tilde{G}_2$: $P$ is a triangle with vertex angles $\frac{\pi}{2}$, $\frac{\pi}{3}$, and $\frac{\pi}{6}$.
\end{enumerate}
These Coxeter systems are  irreducible except for that of type $\tilde{A}_1 \times \tilde{A}_1$, where $W$ is the direct product of two copies of the infinite dihedral group.
\end{example}

Table~\ref{table:wallpaper} gives the correspondences between the  types from Example~\ref{eg:dimension2}, common notations for wallpaper groups, and the quilting piece(s) of shape the polygon $P$. In quilting, both a half-square triangle and a quarter-square triangle are right-angled isosceles triangles, and so correspond to type $\tilde{B}_2/\tilde{C}_2$, but these triangles have different physical properties due either one or two sides being cut on the bias (see, for example, the Glossary to~\cite{Fielke}).

\begin{table}
    \centering
    \renewcommand{\arraystretch}{1.5}
    \begin{tabular}{cccc}
        \bf{Type} & \bf{Conway notation} & \bf{IUC notation} & \bf{Quilting piece(s)} \\ \hline
         $\tilde{A}_1 \times \tilde{A}_1$ & $*2222$ & pmm & Square\\ \hline
        $\tilde{A}_2$ & $*333$ & p3m1 & Equilateral triangle \\ \hline
        $\tilde{B}_2/\tilde{C}_2$ & $*442$ & p4m & Half-square triangle \\
        & & & Quarter-square triangle \\  \hline

        $\tilde{G}_2$ & $*642$ & p6m & 30 degree triangle \\ \hline
    \end{tabular}

    \medskip
    \caption{The types of the $2$-dimensional affine Coxeter systems, from Example~\ref{eg:dimension2}, together with Conway's orbifold notation~\cite{Conway} and the crystallographic notation for this wallpaper group, and corresponding quilting piece(s).}
    \label{table:wallpaper}
\end{table}

\begin{figure}[!ht]
    \centering
    \begin{minipage}{.3\textwidth}
        \centering
        \includegraphics[width=\textwidth]{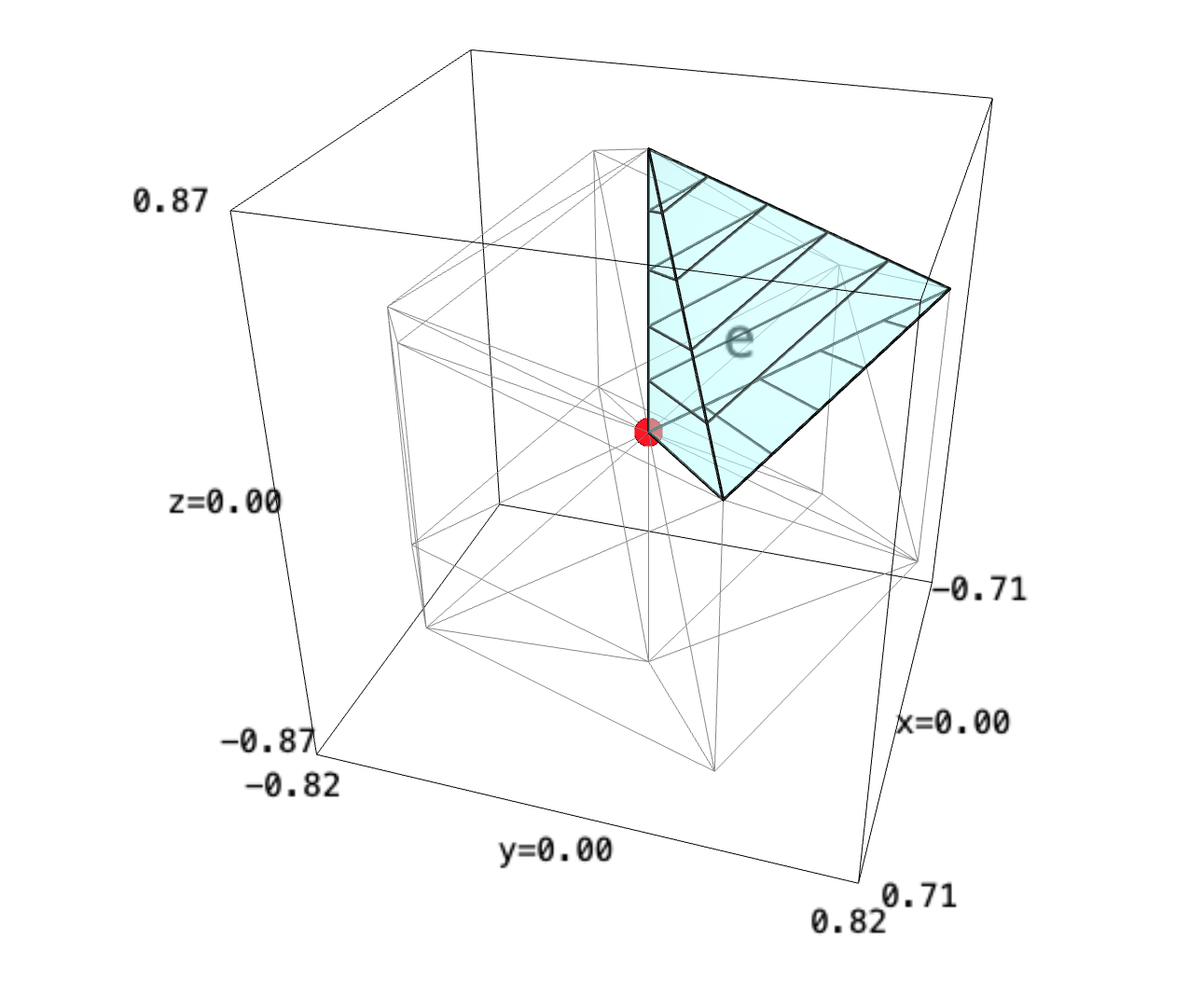}
    \end{minipage}%
    \begin{minipage}{.3\textwidth}
        \centering
        \includegraphics[width=\textwidth]{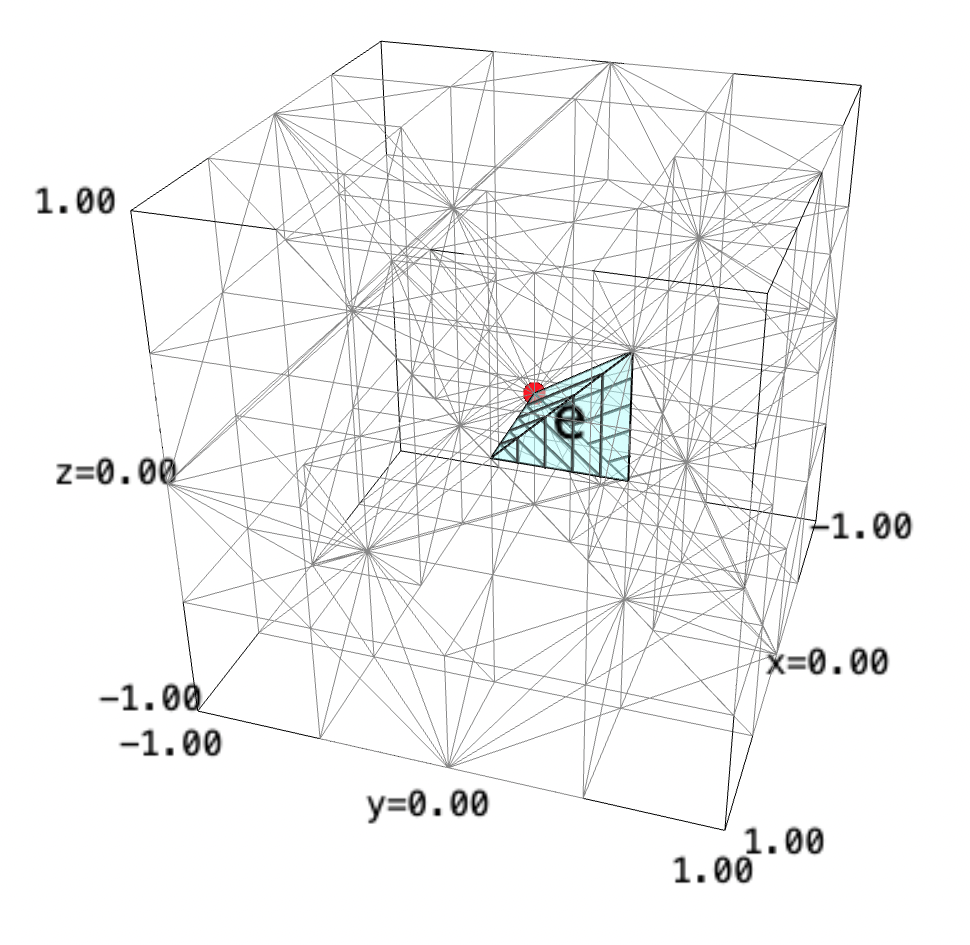}
    \end{minipage}%
    \begin{minipage}{.3\textwidth}
        \centering
        \includegraphics[width=\textwidth]{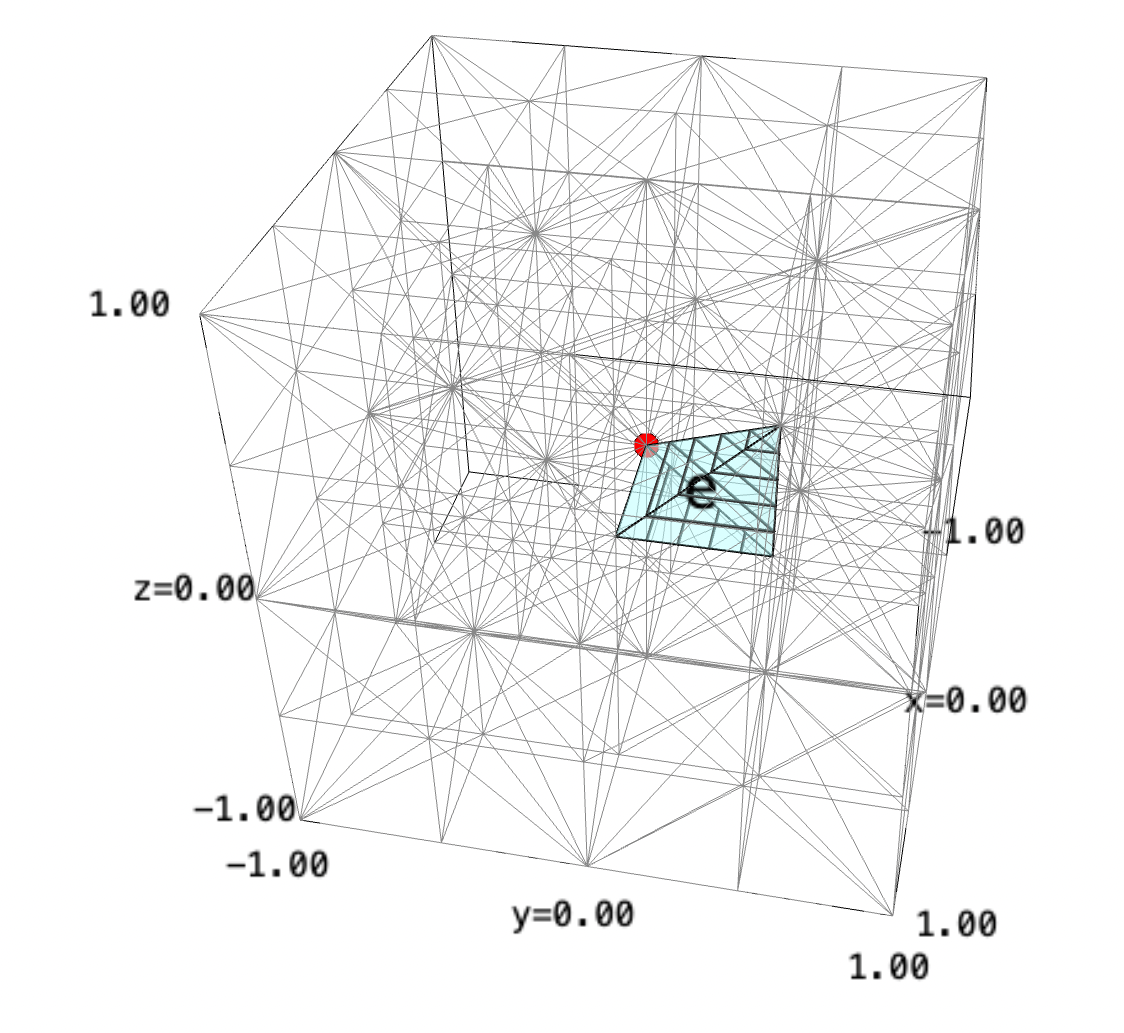}
    \end{minipage}%
    \caption{\footnotesize{From left to right, part of the tesselation of $3$-dimensional Euclidean space by copies of the tetrahedron $P$ in type $\tilde{A}_3$, $\tilde{B}_3$, and $\tilde{C}_3$. The shaded and striped alcove $e$ is the identity element of $W$, and the red dot is the origin.}}
    \label{fig:dimension3_identity}
\end{figure}

\begin{example}\label{eg:dimension3}
The $3$-dimensional irreducible affine Coxeter systems $(W,S)$ are those of types $\tilde{A}_3$, $\tilde{B}_3$ and $\tilde{C}_3$. The corresponding simplex $P$ (a tetrahedron) is different in each of these types. See Figure~\ref{fig:dimension3_identity}.
\end{example}

Let $(W,S)$ be an irreducible affine Coxeter system, with its natural action on $\bE$. As in the introduction, the \emph{coroot lattice} $L$ is the set of all vectors $\lambda \in \bE$ such that the translation by $\lambda$ is in $\aW$. Write $t^\lambda$ for the translation by $\lambda$, so that $\aW$ has translation subgroup $T = \{ t^\lambda \mid \lambda \in L \}$. The coroot lattice has canonical basis elements called the \emph{simple coroots} and denoted $\alpha_1^\vee, \dots, \alpha_n^\vee$, where $n = \dim(\bE)$. For example, the simple coroots $\alpha_1^\vee$ and $\alpha_2^\vee$ are depicted in Figures~\ref{fig:A2tilde_Intro} and~\ref{fig:C2tilde_ConjClass_Intro}. 

Write $\sW$ for the stabilizer of the origin in $\aW$, and define $S_0 = \sW \cap S$. Then $(\sW, S_0)$ is a spherical Coxeter system, and $\sW$ is called the \emph{finite Weyl group}. The alcoves containing the origin are exactly those corresponding to the elements of $\sW$, and every element $x\in\aW$ can be expressed uniquely in the form $x=t^\lambda w$ where $\lambda \in L$ and $w \in \sW$. We then call $t^\lambda$ the \emph{translation part} of $x$, and $w \in \sW$ its \emph{spherical direction}. In fact, $W$ is the semidirect product of its translation subgroup $\{ t^\lambda \mid \lambda \in L \}$ and its finite Weyl group $\sW$, with the action of $\sW$ on translations given by $w t^\lambda w^{-1} = t^{w\lambda}$. If $(W,S)$ is irreducible affine of type~$\tilde{X}_n$, with $S = \{ s_0,s_1, \dots, s_n\}$, then $(\sW, S_0)$ is irreducible spherical of type $X_n$, with $S_0 = \{ s_1,\dots,s_n\}$.

\begin{example}\label{eg:dimension2_spherical}
The finite Weyl groups $W_0$ corresponding to the $2$-dimensional affine Coxeter systems $(W,S)$ from Example~\ref{eg:dimension2} above are as follows.
\begin{enumerate}
    \item Type $A_2$: the group $W_0 = \langle s_1, s_2 \mid s_1^2 = s_2^2 = (s_1 s_2)^3 \rangle$ is the dihedral group of order~$6$, and can also be identified with $\Sym(3)$, the symmetric group on $3$ letters.
    \item Type $B_2$ and type $C_2$: the group $W_0 = \langle s_1, s_2 \mid s_1^2 = s_2^2 = (s_1 s_2)^4 \rangle$ is the dihedral group of order~$8$.
    \item Type $G_2$: the group $W_0 = \langle s_1, s_2 \mid s_1^2 = s_2^2 = (s_1 s_2)^6 \rangle$ is the dihedral group of order~$12$.
    \item Type $A_1 \times A_1$: the group $W_0 = \langle s_1, s_2 \mid s_1^2 = s_2^2 = (s_1 s_2)^2 \rangle$ is the Klein $4$-group i.e. the direct product of two groups of order $2$.
\end{enumerate}
\end{example}

\begin{example}\label{eg:dimension3_spherical}
The finite Weyl group of type $A_3$ has $24$ elements, and can be identified with the symmetric group $\Sym(4)$. The finite Weyl groups in type $B_3$ and type $C_3$ are isomorphic groups with $48$ elements. 
\end{example}

\subsection{The geometry of conjugacy classes and coconjugation sets}\label{sec:MST}

We now briefly recall the relevant results of Mili\'cevi\'c, Schwer and the second author~\cite{MST-affine,MST-Euclidean}. See the introduction for the definitions and our notation for conjugacy classes and coconjugation sets, together with move-sets, fix-sets, and mod-sets.

 Let $(W,S)$ be an $n$-dimensional affine Coxeter system, with coroot lattice $L$ and finite Weyl group $\sW$. We first recall some relationships between the mod-set $\Mod(w)$, the move-set $\Mov(w)$, and the coroot lattice $L$, for elements $w \in \sW$. It is immediate from definitions that $\Mod(w)$ is contained in $\Mov(w)$ and in $L$. Hence $\Mod(w) \subseteq \Mov(w) \cap L$. Moreover, the lattice $L$ can be viewed as a free abelian $\Z$-module (of rank $n$), and then both $\Mod(w)$ and the intersection $\Mov(w) \cap L$ are submodules of $L$. The next statement is a special case of Corollary 1.9(a) of~\cite{MST-affine}.

\begin{thm}\label{thm:ModMove}
    For any $w \in \sW$, the mod-set $\Mod(w)$ has finite index in $\Mov(w) \cap L$.
\end{thm}

\noindent In general, $\Mod(w)$ is a proper submodule of $\Mov(w) \cap L$. For example, in type $\tilde{C}_2$, if $w = s_2$ then $\Mod(w)$ has index $2$ in $\Mov(w) \cap L$ (this explains the ``gaps" between alcoves in the conjugacy class depicted in Figure~\ref{fig:C2tilde_ConjClass_Intro}). 

To fully describe conjugacy classes in $\aW$ we draw on the following statement, which is a special case of Theorem~1.2 of~\cite{MST-Euclidean}.

\begin{thm}[Conjugacy classes]\label{thm:conj}
    Let $x = t^\lambda w \in \aW$, where $\lambda \in L$ and $w \in \sW$. Then the conjugacy class of $x$ in $\aW$ satisfies
    \[
    \xconj = \bigcup_{u \in \sW} t^{u(\lambda + \Mod(w))}uwu^{-1}
    \subseteq \bigcup_{u \in \sW} t^{u(\lambda + \Mov(w))}uwu^{-1}
    =   \bigcup_{u \in \sW} t^{u\Mov(x)}uwu^{-1}.
    \]
\end{thm}

\noindent That is, for any $x = t^\lambda w \in \aW$ the conjugacy class $\xconj$ is found by, for each $u \in \sW$, translating $uwu^{-1}$ by all elements of $u(\lambda + \Mod(w))$. Since $\Mod(w) \subseteq \Mov(w)$, we get the stated containment, and as  $\Mov(x) = \lambda + \Mov(w)$ (by~\cite[Proposition 1.21]{LMPS}), we obtain the final equality.

We next describe a natural decomposition of conjugacy classes. 
In Theorem~\ref{thm:conj}, each subset of the form $t^{u(\lambda + \Mod(w))}uwu^{-1}$, where $u \in \sW$, is called a \emph{component} of the conjugacy class $\xconj$. Recall that $\aW = T \rtimes \sW$, where $T = \{ t^\lambda \mid \lambda \in L \}$ is the translation subgroup of $\aW$. The action of $\aW$ on itself by conjugation induces symmetries between and within components, as follows:

\begin{thm}[Components]\label{thm:components}
    Let $x = t^\lambda w \in \aW$, where $\lambda \in L$ and $w \in \sW$. Then:
    \begin{enumerate}
        \item the conjugation action of $\sW$ permutes the components of $\xconj$; and
        \item the conjugation action of  $T$ induces a transitive action by translations on the elements of each component of $\xconj$.
    \end{enumerate}
\end{thm}

\noindent Part (1) here is immediate from the definition of components, and part (2) is a special case of~\cite[Theorem 2.11(1)]{MST-Euclidean}.

In order to describe coconjugation sets in $\aW$, we need some additional definitions. For $w, w' \in \sW$, the associated \emph{spherical coconjugation set} is given by 
\[
\coconjW w {w'} = \{ u \in \sW \mid uwu^{-1} = w'\} = \coconj{w}{w'} \cap \sW.
\]
Now let $x = t^\lambda w$ and $x' = t^{\lambda'} w'$ be elements of~$\aW$, where $\lambda, \lambda' \in L$ and $w,w' \in \sW$. The corresponding \emph{translation-compatible part} of  $\coconjW{w}{w'}$ is the set
\[
\tcCoconj_{\sW}^{\lambda,\lambda'}(w,w') = \{ u \in \coconjW{w}{w'} \mid \lambda' - u\lambda \in \Mod(w') \}.
\]

The next theorem describes coconjugation sets, and is a special case of Theorem 1.13 of~\cite{MST-Euclidean}.

\begin{thm}[Coconjugation sets]
	\label{thm:coconj} Let $x=t^{\lambda}w$ and $x' = t^{\lambda'}w'$ be elements of $\aW$, where $\lambda, \lambda' \in L$ and $w, w' \in \sW$.  Then 
	\begin{equation}\label{eq:coconjNonempty}
		\coconj{x}{x'} \neq \emptyset \;\; \Longleftrightarrow\;  \tcCoconj_{\sW}^{\lambda,\lambda'}(w,w') \neq \emptyset.
	\end{equation}
	Moreover, if these sets are nonempty, then
	\begin{equation}\label{eq:coconjFix}
	\coconj{x}{x'} = \bigsqcup_{u \in \tcCoconj_{\sW}^{\lambda,\lambda'}(w,w')} t^{\eta_{u}+ (\Fix(w') \cap L)} u 
	\end{equation}
	where for each $u$, the element $\eta_{u} \in L$ is a solution to the equation 
	\begin{equation}\label{eq:coconj}
		(I - w')\eta = \lambda'- u\lambda.
	\end{equation}
\end{thm}

\noindent That is, the coconjugation set $\coconj x {x'}$ has a closed-form description involving $\tcCoconj_{\sW}^{\lambda,\lambda'}(w,w')$ and hence $\Mod(w')$, and its shape is described by translates of $\Fix(w')$. The reason $\Fix(w')$ appears in this statement is that, given one solution $\eta = \eta_u$ to Equation~\eqref{eq:coconj}, all solutions to this equation are of the form $\eta_u + \eta$ where $\eta$ is a solution to the homogeneous equation  $(I - w')\eta = 0$, equivalently, $\eta \in \Fix(w')$.

\begin{rmk} For all irreducible affine Coxeter systems, further results of~\cite{MST-affine} give type-by-type descriptions of mod-sets, and so determine in particular when $\Mod(w) = \Mov(w) \cap L$. We do not use these descriptions in our visualization, instead computing mod-sets directly, so do not recall these results from~\cite{MST-affine} here.
\end{rmk}

\section{Visualization}\label{sec:visualization}

In this section, we describe our app and explain the choices we have made in order to improve the user's understanding of the underlying mathematics.  We describe how to install and run our app in Section~\ref{subsec:app}. Section~\ref{subsec:features} describes its output, especially its visual features, and we explain our choice of colors for alcoves in more detail in Section~\ref{subsec:color}. Finally, we outline our code in Section~\ref{subsec:code}.

\subsection{Application}\label{subsec:app}

In this section we explain how to install and run our app.

The application \texttt{AffineCoxeterExplorer.java}~\cite{app}  is written using Java and SageMath.  Both need to be installed locally before running the app.  Java 11+ can be downloaded from \url{https://adoptium.net}  and SageMath can be downloaded from \url{https://www.sagemath.org} (both are free).  Note that installing SageMath locally can sometimes be difficult.  We found that artificial intelligence can help resolve individual computer issues.  

Once Java and SageMath are installed, the user should download the following files from \url{https://github.com/herron-amy/AffineCoxeterExplorer}:
\begin{enumerate}
    \item \texttt{AffineCoxeterExplorer.java}
    \item \texttt{compute\_helper.sage}
    \item \texttt{README.txt}
    \item one of the following, as appropriate:
        \begin{enumerate}
            \item \texttt{Launch (Linux).sh}
            \item \texttt{Launch (macOS).command}
            \item \texttt{Launch (Windows).bat}
        \end{enumerate}
\end{enumerate}

The user can then run the app by either double clicking on the appropriate Launch option or by inputting 
\begin{center}
     \texttt{java --source 11 AffineCoxeterExplorer.java} 
\end{center}
 in a terminal window, after changing the directory to the one containing the app.  The app will open and the user then chooses the following:
\begin{itemize}
    \item  The type of the affine Coxeter system $(W,S)$ they want to visualize. All $2$-dimensional and all irreducible $3$-dimensional types are available (see Examples~\ref{eg:dimension2} and \ref{eg:dimension3}, respectively).
    \item Whether they want to visualize a conjugacy class or a coconjugation set.
    \begin{itemize}
        \item If ``Conjugacy Class" is selected, the user then inputs one element of $W$.
        \item If ``Coconjugation Set" is selected, the user then inputs two elements of $W$.
    \end{itemize}
    The methods for inputting elements of $W$ are described below.
    \item Size of bounding box, in the range 1--15.  The larger the bounding box, the more alcoves are shown, at the cost of increased run time.  When the user selects the type of $(W,S)$, we suggest a bounding box size which gives a reasonable trade-off, but this size can be modified by the user. 
\end{itemize}
Instead of these selections, the user can instead choose one of 15 quick examples. These illustrate both conjugacy classes and coconjugation sets. 

There are two ways for the user to input elements of the affine Coxeter group $W$, which correspond to two useful algebraic  structures on $W$. The first method corresponds to the presentation of $W$ as a group generated by $S$. The elements of the generating set $S$ for $W$ are indexed by SageMath as $s_0,s_1,\dots, s_n$, with $n = 2$ or $n = 3$ for the cases we consider.  Now given $x = s_{i_1}\dots s_{i_k} = s_{i_1 \dots i_k} \in \aW$, the element~$x$ can be input by entering the string $i_1 \dots i_k$. Note that the word $x = s_{i_1}\dots s_{i_k}$ used to input~$x$ is not required to be reduced. Moreover, even if this word is reduced, it is not required to be the first reduced word in the lexicographic order induced by SageMath's indexing of the generating set $S$. For example, some options for inputting the element $x = s_{0120102}$ from Figure~\ref{fig:A2tilde_Intro} are  \texttt{0120102}, \texttt{0121012} (both of which are reduced), or \texttt{012102121} (which is not reduced). The identity element is input by leaving the entry box empty.

The second method for inputting elements of $W$ corresponds to $W$ being the semidirect product of its translation subgroup by its finite Weyl group $\sW$. For any $x \in \aW$, we have $x = t^\lambda w$ for a unique $\lambda \in L$ and $w \in \sW$. Then $\lambda = c_1 \alpha_1^\vee + \dots + c_n \alpha_n^\vee$ for unique coefficients $c_1,\dots, c_n \in \Z$, where $\alpha_1^\vee,\dots,\alpha_n^\vee$ are the simple coroots. The user can then input the element $x = t^\lambda w$ by entering \texttt{t\_(c\_1,\dots,c\_n)*s\_XX},  where $w = s_{i_1}\dots s_{i_k}$ and $\texttt{XX}$ is the string $i_1 \dots i_k$ (again, this word for $w$ need not be reduced). For example, the element $x = t^{2\alpha_1^\vee + 2\alpha_2^\vee}s_1$ from Figure~\ref{fig:A2tilde_Intro} can be input as: \texttt{t\_(2,2)*s\_1}.

\subsection{Output and visual features}\label{subsec:features}

In this section we describe the output of our app, which comprises both text and graphics. The visual features of the app have been chosen carefully to illustrate the algebraic content of the theorems from Section~\ref{sec:MST} for the general viewer.  For example, the group-theoretic content is represented by color rather than notation, alcoves are labeled directly in the images, and quick pre-coded examples allow the user to view meaningful images without needing to choose group element(s) themselves.  

The text output, which appears in the right-hand column of the app, is the SageMath output. It provides coordinates for the simple coroots, and then a list of all alcoves in the conjugacy class (respectively, coconjugation set) which lie inside the selected bounding box. These alcoves are listed using both their reduced word which is first in lexicographic order, and their translation part (with respect to the basis of simple coroots) and spherical direction. Thus once a conjugacy class $\xconj$ is computed, we have easy access to a list of elements to which $x$ is conjugate.  Any of these elements can then be used in addition to $x$ to visualize a nonempty coconjugation set. 

The graphical output shows the tesselation of the Euclidean plane by alcoves in the $2$-dimensional cases, and by wireframe alcoves (the 1-skeleton of the induced tesselation) in the $3$-dimensional cases. In dimension~$2$ this output is displayed in the center of the app, while in dimension $3$ a browser window will open, and the user can then use the mouse to zoom in or out and rotate the wireframe around, in order to better understand the image. In dimension~$2$, the coroot lattice vertices are indicated by heavy black dots, and the coroot basis vectors $\alpha_1^\vee$ and $\alpha_2^\vee$ are drawn in bold red and bold blue, respectively (see Figures~\ref{fig:A2tilde_Intro} and~\ref{fig:C2tilde_ConjClass_Intro}). We decided not to display this coroot lattice information in dimension $3$, so as to avoid cluttering the image. 

In the graphical output, certain alcoves are decorated to indicate mathematical information. Most importantly, the alcoves in the computed conjugacy class (respectively, coconjugation set) 
that lie within the chosen bounding box are shaded according to their spherical direction. In dimension 2, each spherical direction is assigned its own color; in dimension 3, the color assignment is coarser and differs between conjugacy classes and coconjugation sets, as explained in Section~\ref{subsec:color}.
Shading by spherical direction 
helps the user to see that conjugacy classes lie along images of move-sets, and that coconjugation sets lie along images of fix-sets (as described algebraically by Theorems~\ref{thm:conj} and~\ref{thm:coconj}, respectively). We decided not to depict the move-set or fix-set itself, because the pictures became too cluttered, especially in the 3-dimensional cases, and by Theorem~\ref{thm:ModMove} the move-set can be inferred from the mod-set. 

Our color choices also help the user to distinguish components of conjugacy classes, and to visualise the symmetries between and within components which are described algebraically by Theorem~\ref{thm:components}. For example, on the left of Figure~\ref{fig:A2tilde_Intro} and in the left and central panels of Figure~\ref{fig:A3tilde_Intro}, each component gets assigned a distinct color, and it can be seen that the components may be permuted and that the elements within each component are related to each other by translations. However the elements in distinct components can have the same spherical direction, as in Figure~\ref{fig:C2tilde_ConjClass_Intro}, where both of the horizontal (respectively, vertical) lines of tiles get assigned the same color. In this situation, the spacing between the images of the mod-set corresponding to distinct components serves to distinguish them (if they had the same spherical direction and lay along the same image of the mod-set, they would be the same component).

To help locate the computed conjugacy class (respectively, coconjugation set), in both dimensions $2$ and $3$, the identity alcove is shaded and labeled $e$. If $e$ is an element of the given conjugacy class (respectively, coconjugation set), then $e$ is also striped. This can be seen in Figure~\ref{fig:dimension3_identity}, which depicts the conjugacy class $[e] = \{ e\}$ in the $3$-dimensional cases. In dimension~$2$, we additionally shade  and label all elements of the finite Weyl group $\sW$, and stripe them if they are in the computed conjugacy class (respectively, coconjugation set). In dimension~$3$, since the finite Weyl group $\sW$ has either 24 or 48 elements (see Example~\ref{eg:dimension3_spherical}), such shading and labeling of the entire group $\sW$ is not feasible: there are too many elements of $\sW$ for the viewer to easily distinguish $|\sW|$-many colors, and such labels would overlap each other on a $2$-dimensional screen, so as to become unreadable. Instead, we  indicate the origin in dimension~$3$ by a red dot. 

So that the user can easily compare input and output, the alcove(s) for the user-supplied element(s) are labeled by the corresponding reduced word $s_{i_1 \dots s_{i_k}}$ which is first in lexicographic order, and shaded according to their spherical direction. In addition, for the conjugacy class $\xconj$ and the coconjugation set $\coconj x y$, the alcove $x$ is outlined in bold red, and for the coconjugation set $\coconj x y $, the alcove $y$ is outlined in bold blue.

\subsection{Color choices}\label{subsec:color}

As explained in Section~\ref{subsec:features}, certain alcoves are shaded in our graphical output, according to their spherical direction. In this section, we explain these color choices and what they display to the viewer in more detail.

The $2$-dimensional affine Coxeter groups of types $\tilde A_1\times \tilde A_1$, $\tilde A_2$, $\tilde B_2$, $\tilde C_2$, and $\tilde G_2$ have finite Weyl groups of orders 4, 6, 8, 8, and 12, respectively (see Example~\ref{eg:dimension2_spherical}).  
For these finite Weyl groups, we assign colors chosen for maximum contrast to each element of $\sW$, so that each spherical direction has its own color. 
This assignment is fixed for each type, so that colors may be compared between figures: for example, the purple tiles on the left and right of Figure~\ref{fig:A2tilde_Intro} have the same spherical direction.

In dimension $3$, the finite Weyl groups corresponding to types $\tilde A_3$, $\tilde B_3$, and $\tilde C_3$ contain 24, 48, and 48 elements, respectively (see Example~\ref{eg:dimension3_spherical}). As discussed above, these groups are so big that it is not useful to the viewer to assign a different color to each finite Weyl group element; thus, our color assignment must be coarser. 
To manage this, we observe that if $x = t^\lambda w$ and $x' = t^{\lambda'}w'$ are conjugate elements of $W$, with $\lambda, \lambda' \in L$ and $w,w' \in \sW$, then $w$ and $w'$ must be conjugate in~$\sW$.  It therefore suffices to organize our colors by the conjugacy classes of $W_0$, and we do so in two dual ways, depending on the computation being displayed.

For a conjugacy class, all spherical directions belong to a single conjugacy class of $W_0$, by the observation above. We therefore assign maximally contrasting colors to the elements \emph{within} each conjugacy class of $W_0$, reusing our palette across different classes; thus in a conjugacy class figure, every spherical direction which appears is visually distinguishable. For a coconjugation set, the situation is reversed.  By Theorem~\ref{thm:coconj}, the spherical directions that occur form a coset of the centralizer of $w$ in $W_0$, and such a coset typically meets several conjugacy classes of $W_0$. Here, we instead assign one color to \emph{each} conjugacy class of $W_0$, so that in a coconjugation set figure, alcoves of equal color have conjugate spherical directions. Each of these two schemes spends the available visual contrast on exactly the distinction that the other must sacrifice.  With $W_0$ having 24 or 48 elements, no single coloring scheme can preserve both.

To compute the conjugacy classes of $W_0$, we use injective homomorphisms into symmetric groups. For all $n \geq 1$, the finite Coxeter group of type $A_n$ is isomorphic to the symmetric group $\mathrm{Sym}(n+1)$. Thus in particular, if $(W,S)$ is of type $\tilde{A}_3$, then the finite Weyl group $W_0$ is isomorphic to $\mathrm{Sym}(4)$, where two elements are conjugate if and only if they have the same cycle type. For $W_0$ of type $B_3$ or $C_3$, let $t_1, \ldots, t_5$ be the Coxeter generators for the finite Weyl group of type $A_5$, which is isomorphic to $\mathrm{Sym}(6)$. Then there is an injective homomorphism $\phi \colon W_0 \to \mathrm{Sym}(6)$ induced by $s_1 \mapsto t_1t_5$, $s_2 \mapsto t_2t_4$, and $s_3 \mapsto t_3$.

Regarding types $B_3$ and $C_3$, some care is required because conjugacy in the subgroup $\phi(W_0)$ is finer than conjugacy in $\mathrm{Sym}(6)$: two elements of $\phi(W_0)$ with the same cycle type are conjugate by some element of $\mathrm{Sym}(6)$, but the conjugating permutation need not lie in $\phi(W_0)$. Indeed, the group $W_0$ has ten conjugacy classes, while its image under $\phi$ realizes only eight cycle types. The elements of cycle type $(2,2,1,1)$ form two conjugacy classes of $W_0$, of sizes 6 and 3: geometrically, the first consists of reflections, which fix a plane, and the second consists of rotations by $\pi$ about a coordinate axis, which fix only a line, so no conjugation by isometries can relate them, even though no cycle type computation can tell them apart. Similarly, the elements of cycle type $(2,2,2)$ form two conjugacy classes, of sizes 6 and 1, the singleton being the central element $-I$, which commutes with all of $W_0$ and hence is conjugate only to itself. We therefore use cycle type under these embeddings to separate the eight conjugacy classes that it distinguishes, and separate the remaining two pairs of classes directly.

%Now it is well known that two permutations are conjugate in a symmetric group if and only if they have the same cycle type. Thus in order to determine conjugacy classes in $\sW$ in dimension $3$, it makes sense to find an injective homomorphism from $\sW$ to a symmetric group.   

\subsection{Outline of code}\label{subsec:code}

We conclude by giving an overview of the code for our app.

The code for the $2$-dimensional affine Coxeter groups differs substantially from that for the $3$-dimensional cases.  This is in part because built-in graphics are available for the $2$-dimensional tesselations, but not the $3$-dimensional ones. Also, we display coroot lattice information in dimension $2$ but not in dimension $3$, and our method of assigning colors to spherical directions differs between these dimensions, and in dimension 3 between the two computations, as explained in Section~\ref{subsec:color}. We are, however, able to obtain reduced words and the translation parts and spherical directions for elements of conjugacy classes and coconjugation sets using existing SageMath code for extended affine Weyl groups, in both dimensions $2$ and $3$. In both dimensions, we build labels for alcoves by hand, in order to modify the label font depending on the size of the individual alcoves. 

For the $2$-dimensional cases, we first find Euclidean coordinates for the simple coroots $\alpha_1^\vee$ and $\alpha_2^\vee$ in each type. We then use this basis for $L$ to find all the coroot lattice vertices which are contained within the user-specified bounding box.   

To find the elements of the conjugacy class of $x = t^\lambda w \in \aW$, we start by computing $(w - I)\alpha_i^\vee$ for each simple coroot $\alpha_i^\vee$. The mod-set $\Mod(w) = (w - I)L$ will be the $\Z$-span of these vectors. We then apply each $u \in \sW$ to the user-supplied vector $\lambda$ and to each
$(w - I)\alpha_i^\vee$, and restrict the resulting values of $t^{u\lambda + u\Mod(w)}$ to those within the user-specified bounding box.  
Then for all $u \in \sW$, we conjugate the spherical direction $w$ of the user-supplied element by $u$ as well, to obtain $uwu^{-1}$. Pairing these, we obtain  
 \[ \displaystyle \xconj =[t^\lambda w]=\bigcup_{u\in W_0} t^{u(\lambda + \Mod(w))}u w u^{-1},\]
as in the statement of Theorem~\ref{thm:conj}. 

To find the elements of the coconjugation set $\coconj{x}{x'}$ in dimension $2$, where $x=t^\lambda w$ and $x'=t^{\lambda '} w'$, we first compute the spherical coconjugation set $\coconjW w {w'}$, and then determine its translation-compatible part $\tcCoconj_{\sW}^{\lambda,\lambda'}(w,w')$.  Assuming that $\tcCoconj_{\sW}^{\lambda,\lambda'}(w,w')$ is nonempty, for each $u$ in $\tcCoconj_{\sW}^{\lambda,\lambda'}(w,w')$ we find all solutions $\eta_u$ to the equation $(I - w')\eta = \lambda' - u\lambda$ which lie within our bounding box.
 These $\eta_u$ become the translation part and the $u$ becomes the spherical direction of the elements that are in the coconjugation set $\coconj{x}{x'}$, as given by Theorem~\ref{thm:coconj}.

In $3$ dimensions, we need to build the tesselations for our graphical output by hand. For this, we get the coordinates for the identity alcove from SageMath, then act by the finite Weyl group $\sW$ on this alcove to obtain all alcoves containing the origin.  We now take these finite Weyl group alcoves and translate them across $3$-dimensional space, using the coroot lattice $L$.  The balance of the logic remains the same for the 3-dimensional cases as it does for the 2-dimensional cases.  

In both dimensions $2$ and $3$, we used Claude (version Opus 4.6) by Anthropic, an artificial intelligence source, to provide basic programming assistance and to consolidate the SageMath code into the Java app.    Claude also tweaked the SageMath code to make the app display better.   We tested the app to ensure that its output matches that produced directly by SageMath and corrected discrepancies found in this testing. Claude provided most of the 3-dimensional graphics code, and we tested this code to ensure that it has the correct output. Claude also wrote the README file.  Finally, in revising this paper, we used Claude to suggest and/or format some of the references used for the connection between our visualization and the visual arts. We verified these references and consulted the cited works directly ourselves.

\bibliographystyle{alpha}
\bibliography{refs} 

\newcommand{\etalchar}[1]{$^{#1}$}
\begin{thebibliography}{HdY{\etalchar{+}}02}

\bibitem[Ber]{Berry}
Chris {B}erry {A}rt.
\newblock \url{https://chrisberryart.com/}.
\newblock Accessed: 2026-08-18.

\bibitem[Bod13]{Bodner}
B.~Lynn Bodner.
\newblock The {P}lanar {C}rystallographic {G}roups {R}epresented at the
  {A}lhambra.
\newblock In Craig~S. Kaplan and Reza Sarhangi, editors, {\em Proceedings of
  Bridges 2013: Mathematics, Music, Art, Architecture, Culture}, pages
  225--232, Enschede, Netherlands, 2013. Tessellations Publishing.

\bibitem[Bou02]{Bourbaki}
Nicolas Bourbaki.
\newblock {\em Lie groups and {L}ie algebras. {C}hapters 4--6}.
\newblock Elements of Mathematics (Berlin). Springer-Verlag, Berlin, 2002.
\newblock Translated from the 1968 French original by Andrew Pressley.

\bibitem[CBGS16]{Conway}
John~Horton Conway, Heidi Burgiel, and Chaim Goodman-Strauss.
\newblock {\em The Symmetries of Things}.
\newblock CRC Press, 2 edition, 2016.

\bibitem[CFS09]{crans2009musical}
Alissa~S Crans, Thomas~M Fiore, and Ramon Satyendra.
\newblock Musical actions of dihedral groups.
\newblock {\em The American Mathematical Monthly}, 116(6):479--495, 2009.

\bibitem[Cox79]{Coxeter}
H.~S.~M. Coxeter.
\newblock The non-{E}uclidean symmetry of {E}scher's picture ``{C}ircle {L}imit
  {III}''.
\newblock {\em Leonardo}, 12:19--25, 1979.

\bibitem[Emp19]{Emparan}
Mar{\'i}a~Antonieta Empar{\'a}n.
\newblock The {P}lanar {C}rystallography {G}roups as an {I}conographic
  {A}nalysis {T}ool in {I}slamic {A}rt.
\newblock In David Swart, Kristof Fenyvesi, and Bruce Torrence, editors, {\em
  Proceedings of Bridges 2019: Mathematics, Art, Music, Architecture,
  Education, Culture}, pages 303--310, Linz, Austria, 2019. Tessellations
  Publishing.

\bibitem[Fie11]{Fielke}
Sarah Fielke.
\newblock {\em Quilting}.
\newblock Murdoch Books, Sydney, Australia, 2011.

\bibitem[HdY{\etalchar{+}}02]{Henshilwood}
C.S. Henshilwood, F.~d'Errico, R.~Yates, Z.~Jacobs, C.~Tribolo, G.A.T. Duller,
  N.~Mercier, J.C. Sealy, H.~Valladas, I.~Watts, and A.G. Wintle.
\newblock Emergence of modern human behavior: middle stone age engravings from
  {S}outh {A}frica ({R}eports).
\newblock {\em Science}, 295:1278+, 2002.

\bibitem[Hen]{Henshilwood-photo}
Jegerar og sankarar.
\newblock \url{https://ndla.no/nn/r/historie-vg2/jegere-og-sankere/5a8020a32d}.
\newblock Accessed: 2026-08-19.

\bibitem[Her26]{app}
Amy Herron.
\newblock Affine {C}oxeter {E}xplorer, 2026.
\newblock Version 1.1. Software, archived at
  \url{https://doi.org/10.5281/zenodo.21495268}; source at
  \url{https://github.com/herron-amy/AffineCoxeterExplorer}.

\bibitem[Hum90]{Humphreys}
James~E. Humphreys.
\newblock {\em Reflection groups and {C}oxeter groups}, volume~29 of {\em
  Cambridge Studies in Advanced Mathematics}.
\newblock Cambridge University Press, Cambridge, 1990.

\bibitem[LMPS19]{LMPS}
Joel~Brewster Lewis, Jon McCammond, T.~Kyle Petersen, and Petra Schwer.
\newblock Computing reflection length in an affine {C}oxeter group.
\newblock {\em Trans. Amer. Math. Soc.}, 371(6):4097--4127, 2019.

\bibitem[MST25]{MST-affine}
Elizabeth Mili\'cevi\'c, Petra Schwer, and Anne Thomas.
\newblock The geometry of conjugation in affine {C}oxeter groups.
\newblock {\em Internat. J. Algebra Comput.}, 35(3):403--465, 2025.

\bibitem[MST26]{MST-Euclidean}
Elizabeth Mili\'cevi\'c, Petra Schwer, and Anne Thomas.
\newblock The geometry of conjugation in {E}uclidean isometry groups.
\newblock {\em Enseign. Math.}, 72(3-4):447--466, 2026.

\bibitem[Ron09]{Ronan}
Mark Ronan.
\newblock {\em Lectures on buildings}.
\newblock University of Chicago Press, Chicago, IL, 2009.
\newblock Updated and revised.

\bibitem[Sha23]{Shaw}
Wendy M.~K. Shaw.
\newblock {\em Language Without a Code: Islamic Geometry and Modernity}, pages
  25--48.
\newblock Springer International Publishing, Cham, 2023.

\end{thebibliography}

\end{document}